\documentclass[12pt]{article}

\usepackage[utf8]{inputenc}
\usepackage{t1enc}

\usepackage{amsmath, amsthm, amssymb, amsfonts, amscd}
\usepackage{mathtools}
\usepackage{nccmath}
\usepackage{mathrsfs}
\usepackage[mathscr]{eucal}
\usepackage{stmaryrd}
\usepackage{slashed}
\usepackage{thmtools}
\usepackage{accents}

\usepackage{graphicx, graphics, pict2e, epic, epstopdf}
\usepackage{xcolor}
\usepackage{tikz}
\usetikzlibrary{arrows.meta, positioning}
\usepackage{tikz-cd}
\usepackage{tkz-graph}

\usepackage{geometry}
\usepackage{setspace}
\usepackage{indentfirst}
\usepackage{booktabs}
\usepackage{float}
\usepackage{placeins}
\usepackage{subcaption}
\usepackage{enumerate}
\usepackage[symbol]{footmisc}
\usepackage{lastpage}
\usepackage{authblk}
\usepackage{verbatim}
\usepackage{lipsum}

\usepackage{cite}
\RequirePackage{doi}
\allowdisplaybreaks

\theoremstyle{plain}

\newtheorem{theorem}{Theorem}[section]

\newtheorem{remark}{Remark}

\newtheorem{proposition}{Proposition}

\newtheorem{corollary}{Corollary}

\theoremstyle{remark}

\theoremstyle{definition}
\newtheorem{definition}{Definition}

\newcommand{\E}{\mathcal{E}}
\newcommand{\A}{A}
\newcommand{\R}{\mathbb{R}}
\newcommand{\bmone}{\mathbf{1}}
\newcommand{\Sm}{\mathcal{S}}
\newcommand{\Dm}{\mathcal{D}}
\newcommand{\Ssh}{\mathrm{Ssh}}
\newcommand{\Hck}{\mathcal{H}}
\newcommand{\tr}{\operatorname{tr}}

\newcommand{\kron}{\otimes}

\usepackage{graphicx}

\makeatletter
\renewcommand{\maketitle}{
\begin{center}

{\Large\bfseries \@title\par}
\vspace{6mm}

{\large\bfseries \@author\par}
\vspace{4mm}

{\itshape \@address\par}
\vspace{2mm}

{\small\ttfamily \@email\par}
\vspace{5mm}

\vspace{5mm}

\end{center}
}
\makeatother

\makeatletter
\newcommand{\address}[1]{\gdef\@address{#1}}
\newcommand{\email}[1]{\gdef\@email{#1}}
\makeatother

\address{}
\email{}

\usepackage{fancyhdr}

\usepackage[margin=1cm,%
			font=footnotesize,%
			format=hang,%
			labelsep=period,%
			labelfont=bf]{caption}
\allowdisplaybreaks

\title{
Multiplicativity of Vertex Energy Under the Kronecker Product}
\author{Ronak B. Dudhat$^{a}\footnote{Corresponding author.}$, Vinodray J. Kaneria$^b$, Kalpesh M. Popat$^c$}
\address{$^{a,b,c}$Department of Mathematics, Saurashtra University,\\  Rajkot-360005, Gujarat, India.}
\email{ronakdudhat108@gmail.com, kaneriavinodray@gmail.com,  kalpeshmpopat@gmail.com}

\date{\today}

\begin{document}
\maketitle
\thispagestyle{empty}

\begin{abstract}
We prove that vertex energy is multiplicative under the Kronecker product. Using this core property, we introduce a unified theorem to compute vertex energies for any graph operation obtained through the Kronecker product. As direct applications, we derive exact vertex energy formulas for semishadow, generalized splitting, and shadow-splitting graphs. Furthermore, we completely classify all blueprint matrices that preserve vertex energies.

\vspace{2mm}
\vspace{1mm}

\end{abstract}

\onehalfspacing

\section{Introduction}

 The adjacency matrix of a simple graph $G$ with vertices $v_1,v_2,\dots,v_n$ is the $n\times n$ matrix $A(G)=[a_{ij}]$, where $a_{ij}=1$ if $v_i$ is adjacent to $v_j$, and $0$ otherwise. Let $\lambda_1, \lambda_2, \dots, \lambda_n$ be the eigenvalues of $A(G)$. The energy of $G$, denoted by $\mathcal{E}(G)$, is defined as:
\[
\mathcal{E}(G)=\sum_{i=1}^n |\lambda_i|.
\]

The concept of graph energy was introduced by Ivan Gutman \cite{gutman1978} in 1978. Originally motivated by theoretical chemistry, it provides a mathematical abstraction to approximate the total $\pi$-electron energy of conjugated hydrocarbon molecules within Hückel molecular orbital theory. Over the past decades, this concept has attracted considerable attention and has been studied extensively due to its mathematical richness and broad applications across diverse fields \cite{cvetkovic2010, li2012}.

Graph energy captures the total spectral contribution of a graph, but a finer understanding requires examining how this energy is distributed among individual vertices. This perspective led Arizmendi \emph{et al.} \cite{arizmendi2018} to introduce the concept of \emph{vertex energy} to quantify the spectral contribution of each vertex separately.  For a vertex $v_i$ in a graph $G$, the vertex energy $\mathcal{E}_G(v_i)$ is defined as the $(i,i)$-entry of the matrix absolute value:
$$\mathcal{E}_G(v_i) = |A(G)|_{i,i},$$
where $|A| = (AA^*)^{1/2}$. The total graph energy can be derived by summing the individual vertex energies.
\[\E(G) = \sum_{i=1}^n \E_G(v_i).\]

 Gutman and Furtula \cite{gutmantutorial} showed that the vertex energy can be written in terms of the eigenvalues $\lambda_j$ and the components $C_{j,i}$ of the corresponding orthonormal eigenvectors $C_j$ as:
\vspace{-0.2cm}
$$\E_G(v_i) = \sum_{j=1}^n |\lambda_j| C_{j,i}^2.$$  \vspace{-0.5cm}

Vertex energy has attracted considerable attention in recent years. Existing work includes Coulson-type integral formulas \cite{arizmendi2019}, the behavior of vertex energy in trees \cite{arizmendi2024}, and the classification of small integral graphs \cite{shrikanth2025}. Ramane \emph{et al.} derived the vertex energies of subdivision graphs \cite{ramane2025}.
Balakrishnan~\cite{Balakri2004} proved that the graph energy is multiplicative under kronecker product, i.e. $\mathcal{E}(G\kron H)=\mathcal{E}(G)\, \mathcal{E}(H).$
The Kronecker product gives rise to several graph operations whose adjacency matrices are Kronecker products of $\A(G)$ with small matrices built from $I$, $J$ and $\bmone$.  The energies of the $m$-splitting and $m$-shadow graphs were obtained by Vaidya and Popat~\cite{kmp}, and their vertex energies were later determined by Harshitha \emph{et al.}~\cite{harshithavertex}, establishing that the vertex energy remains invariant under $m$-shadow graphs. Furthermore, Dede and Popat proved that the double graphs and bipartite double covers preserve the vertex energy \cite{dede2025}. The energy of the $p$-semishadow graph was obtained by Zahra \emph{et al.}~\cite{zahra2025}, and that of the $(p,q)$-generalized splitting and $(c,k)$-shadow-splitting graphs by Dudhat \emph{et al.}~\cite{Dudhat2026}. 

In this paper we prove multiplicativity of vertex energy under the kronecker product and as its application we deduce the vertex energy formulas for the $p$-semishadow, $(p,q)$-generalized splitting and $(c,k)$-shadow-splitting graphs, and also characterize vertex energy invariant graph operations derived from kronecker product, which extends the results of~\cite{harshithavertex} and \cite{dede2025}. 

\section{Preliminaries}\label{sec:operations}

Throughout this paper, $G$ denotes a simple, connected and undirected graph with vertex set $V(G)=\{v_1,\dots,v_n\}$. For $v\in V(G)$, $N_G(v)$ is its neighbourhood. Let $I_n$ denote the $n\times n$ identity matrix, $J_{m\times n}$ the $m\times n$ all-ones matrix (written $J_n$ when $m=n$), and $\bmone_n$ the all-ones column vector in $\R^{n}$. For matrices $A=(a_{ij})\in\R^{m\times n}$ and $B\in\R^{p\times q}$, the \emph{Kronecker product} $A\kron B$ is the $mp\times nq$ block matrix whose $(i,j)$ block is $a_{ij}B$, that is,
$A\kron B=(a_{ij}B)_{1\le i\le m,\;1\le j\le n}$. We will use the following standard facts about the Kronecker product
\cite{Horn1991}. If $P,Q,R,S$ are matrices of compatible sizes, then
$(P\kron Q)(R\kron S)=(PR)\kron(QS)$. Moreover, if $A$ and $B$ are symmetric with eigenpairs $(\lambda_i,X_i)$ and $(\mu_j,Y_j)$, where the eigenvectors are orthonormal, then
$A\kron B$ has eigenpairs $(\lambda_i\mu_j,\,X_i\kron Y_j)$, and these eigenvectors form a complete orthonormal set. The direct sum $A \oplus B$ is the block diagonal matrix with diagonal blocks $A$ and $B$; more generally, $\displaystyle\bigoplus_{i=1}^k A_i$ denotes the block diagonal matrix with diagonal blocks $A_1, \dots, A_k$.

\begin{definition}[\cite{Weichsel1962}]
For simple graphs $G$ and $H$, the \emph{Kronecker product} $G\kron H$ is a graph with vertex set $V(G)\times V(H)$ and edge set
$\{(u_1,v_1)(u_2,v_2):u_1u_2\in E(G),\ v_1v_2\in E(H)\}$.
\end{definition}

\begin{definition}[\cite{AbdelAal2013}]
Let $m\ge 1$. The \emph{$m$-splitting graph} $\Sm_m(G)$ is obtained from $G$ by adding $m$ independent sets $U^{(k)}=\{u_1^{(k)},\dots,u_n^{(k)}\}$, $k=1,\dots,m$, where $N(u_i^{(k)})=N_G(v_i)$. The \emph{$m$-shadow graph} $\Dm_m(G)$ consists of $m$ copies of $G$, where $N(u_i^{(k)})=\{u_j^{(\ell)} : v_j\in N_G(v_i),\ \ell=1,\dots,m\}$.
\end{definition}

\begin{definition}[\cite{Patel2023}]
Let $p\ge 1$. The \emph{$p$-semishadow graph} $\Ssh_p(G)$ is formed from $G$ and $p$ copies $G^{(1)},\dots,G^{(p)}$ of $G$, where $u_i^{(r)}$ denotes the copy of $v_i$ in $G^{(r)}$, with
\[
N(u_i^{(r)})=\bigcup_{v_j\in N_G(v_i)}\{u_j^{(r)},\,v_j\} \quad \text{and} \quad
N(v_i)=\bigcup_{v_j\in N_G(v_i)}\{v_j,\,u_j^{(1)},\dots,u_j^{(p)}\}.
\]
\end{definition}

\begin{definition}[\cite{Dudhat2026}]
Let $p,q\ge 1$. The \emph{$(p,q)$-generalized splitting graph} $\Sm_{p,q}(G)$ is formed from $p$ copies $G^{(1)},\dots,G^{(p)}$ of $G$ together with $q$ independent sets $U^{(k)}=\{u_1^{(k)},\dots,u_n^{(k)}\}$, $k=1,\dots,q$, where
\[
N(u_i^{(k)})=\bigcup_{\ell=1}^{p}N_{G^{(\ell)}}(v_i^{(\ell)}).
\]
\end{definition}

\begin{definition}[\cite{Dudhat2026}]
For $c,k\ge 1$, the \emph{$(c,k)$-shadow-splitting graph} $\Hck_{c,k}(G)$ is formed from $c$ copies of $G$ and $k$ independent sets, where the copies are joined to one another by the shadow relation and to the added sets by the splitting relation.
\end{definition}

The adjacency matrices of these operations take a simple form once the vertices are ordered block by block: first the vertices of $G$ (or of its first copy), then those of the remaining copies, and finally those of the added independent sets, each block following the order $v_1,\dots,v_n$. With this ordering, each adjacency matrix is a Kronecker product:
\begin{align*}
\A(G\kron H) &= \A(G) \kron \A(H),
&\A\big(\Sm_m(G)\big) &= \begin{bmatrix} 1 & \bmone_m^{T}\\ \bmone_m & 0\end{bmatrix}\kron \A(G), \\[2mm]
\A\big(\Dm_m(G)\big) &= J_m\kron \A(G),
&\A\big(\Ssh_p(G)\big) &= \begin{bmatrix} 1 & \bmone_p^{T}\\ \bmone_p & I_p\end{bmatrix}\kron \A(G), \\[2mm]
\A\big(\Sm_{p,q}(G)\big) &= \begin{bmatrix} I_p & J_{p\times q}\\ J_{q\times p} & 0\end{bmatrix}\kron \A(G),
&\A\big(\Hck_{c,k}(G)\big) &= \begin{bmatrix} J_c & J_{c\times k}\\ J_{k\times c} & 0\end{bmatrix}\kron \A(G).
\end{align*}

\section{Main Results}
We begin with the vertex-energy of the Kronecker product of two graphs, which motivates the more general framework developed in this paper.
\begin{proposition} \label{prop:kron}
Let $G$ and $H$ be simple graphs. For every vertex $(u,v)\in V(G\kron H)$,
\begin{equation}\label{eq:kronvertex}
\E_{G\kron H}(u,v)=\E_G(u)\,\E_H(v).
\end{equation}
\end{proposition}
\begin{proof}
Let $G$ and $H$ have orders $n_1$ and $n_2$, with orthonormal eigenpairs $(\lambda_j,X_j)$ and
$(\mu_\ell,Y_\ell)$. By Gutman and Furtula,
$\E_G(u_i)=\displaystyle\sum_j|\lambda_j|X_{j,i}^2$ and
$\E_H(v_k)=\displaystyle\sum_\ell|\mu_\ell|Y_{\ell,k}^2$. For the matrix
$\A(G\kron H)=\A(G)\kron\A(H)$ the orthonormal eigenpairs are
$(\lambda_j\mu_\ell,\,X_j\kron Y_\ell)$, and the coordinate of
$X_j\kron Y_\ell$ at $(u_i,v_k)$ equals $X_{j,i}Y_{\ell,k}$. Hence
\[
\E_{G\kron H}(u_i,v_k)
=\sum_{j,\ell}|\lambda_j\mu_\ell|\,(X_{j,i}Y_{\ell,k})^2
=\Big(\sum_j|\lambda_j|X_{j,i}^2\Big)\Big(\sum_\ell|\mu_\ell|Y_{\ell,k}^2\Big)
=\E_G(u_i)\E_H(v_k).
\]
\end{proof}
Summing over all vertices gives  $\E(G\kron H)=\E(G)\E(H),$ which recovers the multiplicativity of the energy established by Balakrishnan \cite{Balakri2004}. The proof never uses the fact that $\A(H)$ is the adjacency matrix of a graph, so $\A(H)$ may be replaced by an arbitrary symmetric matrix. For a symmetric matrix $M\in\R^{N\times N}$ with orthonormal eigenpairs $(\theta_s,y_s)$, we define the \emph{energy of index $t$} by
\begin{equation}\label{eq:Mvertexenergy}
\E_M(t)=|M|_{tt}=\sum_{s=1}^{N}|\theta_s|\,(y_s)_t^{2},
\qquad
\E(M)=\operatorname{tr}|M|=\sum_{t=1}^{N}\E_M(t).
\end{equation}
For $M=\A(H)$ this is the usual vertex energy of $H$.

\begin{theorem}\label{thm:master}
Let $G$  be a simple graph of order $n$ and let $M\in\{0,1\}^{N\times N}$ be symmetric.
Let $B=M\kron\A(G)$, and partition the index set of $B$ into $N$ blocks
$V^{(1)},\dots,V^{(N)}$,  each corresponding to $\{v_1,\dots,v_n\}$.Writing $(t,i)$ for the
$i$-th vertex of block~$t$ we have
\begin{equation}\label{eq:masterformula}
\E_B\big((t,i)\big)=\E_M(t)\,\E_G(v_i),
\qquad t=1,\dots,N,\ \ i=1,\dots,n.
\end{equation}
\end{theorem}
\begin{proof}
Let $(\theta_s,y_s)$ and $(\lambda_j,X_j)$ be orthonormal eigenpairs of $M$ and $\A(G)$, respectively. By \eqref{eq:Mvertexenergy} and Gutman--Furtula formula, $\E_M(t)=\sum_s|\theta_s|(y_s)_t^2$ and $\E_G(v_i)=\sum_j|\lambda_j|X_{j,i}^2$. The matrix $B=M\kron\A(G)$, has the orthonormal eigenpairs $(\theta_s\lambda_j,\,y_s\kron X_j)$, and the coordinate of $y_s\kron X_j$ at $(t,i)$ is $(y_s)_tX_{j,i}$. Therefore
\[
\E_B\big((t,i)\big)
=\sum_{s,j}|\theta_s\lambda_j|\,\big((y_s)_t X_{j,i}\big)^2
=\Big(\sum_s|\theta_s|(y_s)_t^2\Big)\Big(\sum_j|\lambda_j|X_{j,i}^2\Big)
=\E_M(t)\,\E_G(v_i).
\]
\end{proof}
\begin{remark}
Summing \eqref{eq:masterformula} over all $(t,i)$ gives $\E(B)=\E(M)\,\E(G)$. If $M$ is a symmetric $0/1$ matrix, then $B$ is the adjacency matrix of a simple graph, and Proposition~\ref{prop:kron} is the case $M=\A(H)$. Moreover, the factor $\E_M(t)$ depends only on $M$ and not on $G$, so a fixed operation $G\mapsto M\kron\A(G)$ rescales the vertex energies of $G$ by the same constants for every base graph $G$.
\end{remark} 

\subsection{Applications to Graph Operations}\label{sec:cor}
We now apply Theorem~\ref{thm:master} to the operations of Section~\ref{sec:operations}. Each operation is of the form $M\kron\A(G)$ for a suitable symmetric matrix $M$, so it suffices to identify $M$ and compute $\E_M(t)$ from \eqref{eq:Mvertexenergy}.

Recall from Arizmendi et al.~\cite{arizmendi2018} that for any symmetric matrix $M$, the diagonal entries of $|M|$ are given by $|M|_{tt} = \displaystyle\sum_{\lambda} p_t(\lambda)|\lambda|$, where the spectral weights $p_t(\lambda) \ge 0$ satisfy the moment equations
\[\sum_{\lambda} p_t(\lambda)\lambda^k = (M^k)_{tt} \quad\text{for}\quad k = 0, 1, 2.\]

\begin{corollary}\label{cor:shadow}
For $\Dm_m(G),\;\; \E_{\Dm_m(G)}(w)=\E_G(v_i)$, for every $m\ge 1$ and every vertex $w=u_i^{(k)}$.
\end{corollary}

\begin{proof}
For $M=J_m$, eigenvalues are $m \ge 0$ and $0$. Since all eigenvalues are non-negative, $|\lambda| = \lambda$ for all $\lambda$, which immediately yields
\[
|M|_{tt} = \sum_\lambda p_t(\lambda)|\lambda| = \sum_\lambda p_t(\lambda)\lambda = M_{tt} = (J_m)_{tt} = 1 \qquad \text{for all } t = 1, \dots, m.
\]
Theorem~\ref{thm:master} proves $\E_{\Dm_m(G)}(u_i^{(k)}) = |M|_{kk}\E_G(v_i) = \E_G(v_i)$.
\end{proof}

\begin{corollary}\label{cor:psplit}
For $\Sm_p(G)$,
\[
\E_{\Sm_p(G)}(w)=
\begin{cases}
\dfrac{1+2p}{\sqrt{1+4p}}\,\E_G(v_i),
& w=v_i\in V(G),\\[0.5cm]
\dfrac{2}{\sqrt{1+4p}}\,\E_G(v_i),
& w=u_i^{(k)},\ k=1,\dots,p.
\end{cases}
\]
\end{corollary}

\begin{proof}
For \(
\begin{medsize}
M = \begin{pmatrix} 1 & \bmone_p^T \\ \bmone_p & 0_p \end{pmatrix}, \;
M^2 = \begin{pmatrix} 1+p & \bmone_p^T \\ \bmone_p & J_p \end{pmatrix}
\end{medsize}
\), and its non-zero eigenvalues are $\alpha_\pm = \mfrac{1\pm D}{2}$ with $D = \sqrt{1+4p}$, satisfying $\alpha_++\alpha_-=1, \; \alpha_+\alpha_-=-p$, and $|\alpha_\pm| = \pm\alpha_\pm$.

For slot $1$ ($w = v_i$), $p_0 = 0$, and the equations $p_+ + p_- = 1$ and $p_+\alpha_+ + p_-\alpha_- = 1$ immediately give $p_\pm = \pm\mfrac{\alpha_\pm}{D}$. Thus,
\[
|M|_{11} = p_+\alpha_+ - p_-\alpha_- = \mfrac{\alpha_+^2 + \alpha_-^2}{D} = \mfrac{1+2p}{\sqrt{1+4p}}.
\]
For slot $1+k$ ($w = u_i^{(k)}$), $p_+\alpha_+ + p_-\alpha_- = 0$ and $p_+\alpha_+^2 + p_-\alpha_-^2 = (M^2)_{1+k,1+k} = 1$ reduce to $p(p_+ + p_-) = 1$, which gives $p_\pm = \mp\mfrac{\alpha_\mp}{pD}$. Thus,
\[
|M|_{1+k,1+k} = p_+\alpha_+ - p_-\alpha_- = \mfrac{2(-\alpha_+\alpha_-)}{pD} = \mfrac{2}{\sqrt{1+4p}}.
\]
Applying Theorem~\ref{thm:master} completes the proof.
\end{proof}

\begin{remark}
The results of Corollaries 1 and 2 recover the vertex energy formulas previously obtained by Harshitha \emph{et al.}~\cite{harshithavertex}. Our framework demonstrates that both formulas arise as immediate, one-line consequences of Theorem~\ref{thm:master}, offering a much more transparent and streamlined derivation.
\end{remark}

\begin{corollary}\label{cor:ssh}
For $\Ssh_p(G)$,
\[
\E_{\Ssh_p(G)}(w)=
\begin{cases}
\sqrt{p}\,\E_G(v_i), & w=v_i\in V(G),\\[0.5cm]
\left(\mfrac{1}{\sqrt p}+\mfrac{p-1}{p}\right)\E_G(v_i),
& w=u_i^{(r)},\ r=1,\dots,p.
\end{cases}
\]
\end{corollary}

\begin{proof}
For $M=\begin{medsize}\begin{pmatrix}1&\bmone_p^{T}\\ \bmone_p&I_p\end{pmatrix}\end{medsize}, \, M^2=\begin{medsize}\begin{pmatrix}1+p&2\bmone_p^{T}\\ 2\bmone_p&I_p+J_p\end{pmatrix}\end{medsize}$, and its distinct eigenvalues are $\lambda_1 = 1+\sqrt{p}$, $\lambda_2 = 1-\sqrt{p}$, and $\lambda_3 = 1$. For $p \ge 1$, we have $|\lambda_1| = 1+\sqrt{p}$, $|\lambda_2| = \sqrt{p}-1$, and $|\lambda_3| = 1$. The relations $p_1 + p_2 + p_3 = 1$ and $p_1\lambda_1 + p_2\lambda_2 + p_3\lambda_3 = M_{tt} = 1$ give $p_1 = p_2$, so the shifted second moment yields \[
p_1(\lambda_1 - 1)^2 + p_2(\lambda_2 - 1)^2 + p_3(0)^2 = (M^2)_{tt} - 2M_{tt} + 1 \implies 2p\,p_1 = (M^2)_{tt} - 1.
\]

For slot $1$ ($w = v_i$), $(M^2)_{11} = 1+p$, so $2p\,p_1 = p$ gives $p_1 = p_2 = \mfrac{1}{2}$ and $p_3 = 0$. Thus,
\[
|M|_{11} = \mfrac{1}{2}(1+\sqrt{p}) + \mfrac{1}{2}(\sqrt{p}-1) = \sqrt{p}.
\]
For slot $r+1$ ($w = u_i^{(r)}$, $r=1,\dots,p$), $(M^2)_{r+1,r+1} = 2$, so $2p\,p_1 = 1$ gives $p_1 = p_2 = \mfrac{1}{2p}$ and $p_3 = \mfrac{p-1}{p}$. Thus,
\[
|M|_{r+1,r+1} = \mfrac{1}{2p}(1+\sqrt{p}) + \mfrac{1}{2p}(\sqrt{p}-1) + \mfrac{p-1}{p}(1) = \mfrac{1}{\sqrt{p}} + \mfrac{p-1}{p}.
\]
Applying Theorem~\ref{thm:master} completes the proof.
\end{proof}

\begin{remark}
Summation gives $\sqrt{p}\,\E(G) + p\left(\mfrac{1}{\sqrt{p}} + \mfrac{p-1}{p}\right)\E(G) = (2\sqrt{p}+p-1)\E(G)$, which recovers the energy of $\Ssh_p(G)$ obtained by Zahra \emph{et al.}~\cite{zahra2025}.
\end{remark}

\begin{corollary}\label{cor:split}
For $\Sm_{p,q}(G)$,
\[
\E_{\Sm_{p,q}(G)}(w)=
\begin{cases}
\left(\dfrac{1+2pq}{p\sqrt{1+4pq}}+\dfrac{p-1}{p}\right)\E_G(v_i),
& w=v_i^{(\ell)},\ \ell=1,\dots,p,\\[0.5cm]
\dfrac{2p}{\sqrt{1+4pq}}\,\E_G(v_i),
& w=u_i^{(k)},\ k=1,\dots,q.
\end{cases}
\]
\end{corollary}

\begin{proof}
For $M=\begin{medsize}\begin{pmatrix} I_p & J_{p\times q}\\ J_{q\times p} & 0_q\end{pmatrix}, 
\;
M^2 = \begin{pmatrix} I_p + q J_p & J_{p\times q}\\ J_{q\times p} & p J_q\end{pmatrix}
\end{medsize}$, and its distinct eigenvalues are $\lambda=1$, $\lambda=0$, and $\alpha_\pm = \mfrac{1\pm D}{2}$ (with $D = \sqrt{1+4pq}$). They satisfy $\alpha_++\alpha_-=1$, $\alpha_+\alpha_-=-pq$, and $\alpha_\pm^2 - \alpha_\pm = pq$. Since $\alpha_+ > 0 > \alpha_-$, we have $|\alpha_\pm| = \pm\alpha_\pm$.

For slot $\ell$ ($w=v_i^{(\ell)}$, $\ell=1,\dots,p$),\, $p_0 = 0$. Using the first and second moment equations $pq(p_+ + p_-) = (M^2)_{\ell\ell} - M_{\ell\ell} = q$, so $p_+ + p_- = \mfrac{1}{p}$ and $p_1 = \mfrac{p-1}{p}$. Since $p_+(\alpha_+-1) + p_-(\alpha_--1) = 0$, we get $p_+\alpha_- = -p_-\alpha_+$, and $p_+ = \mfrac{\alpha_+}{pD}$ and $p_- = \mfrac{-\alpha_-}{pD}$. Thus,
\[
|M|_{\ell\ell} = p_+\alpha_+ - p_-\alpha_- + p_1(1) = \mfrac{\alpha_+^2 + \alpha_-^2}{pD} + \mfrac{p-1}{p} = \mfrac{1+2pq}{p\sqrt{1+4pq}} + \mfrac{p-1}{p}.
\]
For slot $p+k$ ($w=u_i^{(k)}$, $k=1,\dots,q$), $p_1 = 0$. The first moment gives $p_+\alpha_+ + p_-\alpha_- = 0$, so $(M^2)_{p+k,p+k} = p$ reduces to $pq(p_+ + p_-) = p$, which yields $p_+ = \mfrac{-\alpha_-}{qD}$ and $p_- = \mfrac{\alpha_+}{qD}$. Thus,
\[
|M|_{p+k,p+k} = p_+\alpha_+ - p_-\alpha_- = \mfrac{2(-\alpha_+\alpha_-)}{qD} = \mfrac{2pq}{qD} = \mfrac{2p}{\sqrt{1+4pq}}.
\]
Applying Theorem~\ref{thm:master} completes the proof.
\end{proof}

\begin{corollary}\label{cor:shadowsplit}
For $\Hck_{c,k}(G)$,
\[
\E_{\Hck_{c,k}(G)}(w)=
\begin{cases}
\dfrac{c+2k}{\sqrt{c^2+4ck}}\,\E_G(v_i),
& w=v_i^{(\ell)},\ \ell=1,\dots,c,\\[0.5cm]
\dfrac{2c}{\sqrt{c^2+4ck}}\,\E_G(v_i),
& w=u_i^{(r)},\ r=1,\dots,k.
\end{cases}
\]
\end{corollary}

\begin{proof}
For $M=\begin{medsize}\begin{pmatrix} J_c & J_{c\times k}\\ J_{k\times c} & 0_k\end{pmatrix},
\;
M^2 = \begin{pmatrix} (c+k)J_c & c J_{c\times k}\\ c J_{k\times c} & c J_k\end{pmatrix}
\end{medsize}$, and its distinct eigenvalues are $\lambda=0$ and $\alpha_\pm = \mfrac{c\pm D}{2}$ (with $D = \sqrt{c^2+4ck}$). They satisfy $\alpha_++\alpha_-=c$, $\alpha_+\alpha_-=-ck$, and $\alpha_\pm^2 - c\alpha_\pm = ck$. $|\alpha_\pm| = \pm\alpha_\pm$ as $\alpha_+ > 0 > \alpha_-$.

For slot $\ell$ ($w=v_i^{(\ell)}$, $\ell=1,\dots,c$), by first and second moment $ck(p_+ + p_-) = (M^2)_{\ell\ell} - c M_{\ell\ell} = k$, so $p_+ + p_- = \mfrac{1}{c}$ and $p_0 = \mfrac{c-1}{c}$. Since $p_+(\alpha_+-c) + p_-(\alpha_--c) = 0$, we get $p_+\alpha_- = -p_-\alpha_+$, and $p_+ = \mfrac{\alpha_+}{cD}$ and $p_- = \mfrac{-\alpha_-}{cD}$. Thus,
\[
|M|_{\ell\ell} = p_+\alpha_+ + p_-(-\alpha_-) + p_0(0) = \mfrac{\alpha_+^2 + \alpha_-^2}{cD} = \mfrac{c^2 + 2ck}{cD} = \mfrac{c+2k}{\sqrt{c^2+4ck}}.
\]
For slot $c+r$ ($w=u_i^{(r)}$, $r=1,\dots,k$), the first moment gives $p_+\alpha_+ + p_-\alpha_- = 0$, so $(M^2)_{c+r,c+r} = c$ reduces to $ck(p_+ + p_-) = c$, which yields $p_+ = \mfrac{-\alpha_-}{kD}$ and $p_- = \mfrac{\alpha_+}{kD}$. Thus,
\[
|M|_{c+r,c+r} = p_+\alpha_+ + p_-(-\alpha_-) = \mfrac{2(-\alpha_+\alpha_-)}{kD} = \mfrac{2ck}{kD} = \mfrac{2c}{\sqrt{c^2+4ck}}.
\]
Applying Theorem~\ref{thm:master} completes the proof.
\end{proof}

\begin{remark}
Summation recovers $\E(\Sm_{p,q}(G))=(\sqrt{1+4pq}+p-1)\E(G)$ and $\E(\Hck_{c,k}(G))=\sqrt{c^2+4ck}\,\E(G)$ of Dudhat et al.\ \cite{Dudhat2026}.
\end{remark}

\subsection{Classification of Vertex Energy-Preserving Operations}

In the operations examined in the previous subsection, the vertex energies of the derived graphs scale by distinct factors depending on the slot $t$. On the other hand, as noted in the introduction, certain operations---such as the $m$-shadow graph~\cite{harshithavertex} and bipartite double covers~\cite{dede2025}---leave the vertex energies completely invariant. 

Within our framework, this motivates the general inverse question: \emph{Which matrices leave the vertex energy unconditionally invariant across all vertices in the Kronecker product?}

By Theorem~\ref{thm:master}, for the derived graph with adjacency matrix $B = M \kron \A(G)$ the vertex energy of $G$ remains invariant across all vertices if and only if
\(
|M|_{tt} = 1, \;\; \forall \; t = 1, \dots, N.
\)
We now establish a complete classification of all such graphs, beginning with connected graphs.

\begin{theorem}\label{thm:connected_matrix_classification}
Let $M \in \{0,1\}^{N \times N}$ be a symmetric matrix whose underlying graph is connected. Then $|M|_{tt} = 1$ for all $t = 1, \dots, N$ if and only if
\[
M \cong \A(K_{r,r}) \quad (r \ge 1, \; N = 2r) \qquad \text{or} \qquad M \cong J_N \quad (N \ge 1).
\]
\end{theorem}

\begin{proof}
 If $M \cong \A(K_{r,r})$, then $M^2 = r(J_r \oplus J_r)$, so its unique positive semidefinite square root is $|M| = J_r \oplus J_r$, giving $|M|_{tt} = (J_r)_{tt} = 1$. If $M \cong J_N$, since $J_N$ is positive semidefinite, $|J_N| = J_N$, which directly gives $|M|_{tt} = (J_N)_{tt} = 1$ for all $t$.

 Conversely, let $M \in \{0,1\}^{N \times N}$ be symmetric matrix without zero diagonal block and $\E_M(t) = |M|_{tt} = 1$ for all $t = 1, \dots, N$. Partition the vertex set into $V_0 = \{t : M_{tt} = 0\}$ and $V_1 = \{t : M_{tt} = 1\}$. For any $v \in V_1$, we have $M_{vv} = 1$ and $|M|_{vv} = 1$, which gives
\[
|M|_{vv} - M_{vv} = \sum_{\lambda_i < 0} 2|\lambda_i|\, u_{i,v}^2 = 1 - 1 = 0,
\]
forcing $u_{i,v} = 0$ for every negative eigenvalue $\lambda_i < 0$. Hence all negative eigenvectors vanish on $V_1$, and equating blocks in $|M|^2 = M^2$ leaves no edges between $V_0$ and $V_1$; by connectivity, the diagonal entries must be uniformly zero ($V_1 = \emptyset$) or uniformly one ($V_0 = \emptyset$).

If $M_{tt} = 0$ for all $t$, then $M$ represents a connected simple graph. Let $\lambda_1 \ge \dots \ge \lambda_N$ denote the eigenvalues of $M$ with corresponding orthonormal eigenvectors $u_1, \dots, u_N$, and let $\lambda_{\max} = \displaystyle\max_{1 \le i \le N} |\lambda_i|$. Since $|\lambda_i|(\lambda_{\max} - |\lambda_i|) \ge 0$, weighting by $u_{i,t}^2 \ge 0$ and summing over all $i = 1, \dots, N$ yields

\[
\lambda_{\max} \E_M(t) - (M^2)_{tt} = \sum_{i=1}^N |\lambda_i|(\lambda_{\max} - |\lambda_i|)\, u_{i,t}^2 \ge 0.
\]
Because $\E_M(t) = 1$ and $(M^2)_{tt} = \deg(t)$, we obtain $\deg(t) \le \lambda_{\max}$ for every vertex $t$. Also, by combining $\E(M) = N$ with $\tr(M) = 0$, we get
\[
\sum_{\lambda_i > 0} \lambda_i = \sum_{\lambda_i < 0} |\lambda_i| = \mfrac{1}{2}\E(M) = \mfrac{N}{2}.
\]
Because $\lambda_{\max} \le \displaystyle\sum_{\lambda_i > 0} \lambda_i = \mfrac{N}{2}$, we have $\deg(t) \le \lambda_{\max} \le \mfrac{N}{2}$ for all $t$. Since each eigenvalue satisfies $|\lambda_i| \le \lambda_{\max} \le \mfrac{N}{2}$, their sums can equal $\mfrac{N}{2}$ if and only if $\lambda_1 = -\lambda_N = \mfrac{N}{2}$ and all intermediate eigenvalues vanish i.e. $\lambda_2 = \dots = \lambda_{N-1} = 0$.
Setting $r = \mfrac{N}{2}$, the non-zero spectrum $\{r^1, (-r)^1\}$ and regularity $\deg(t) = \lambda_{\max} = r$ uniquely characterize $M$ as $M \cong \A(K_{r,r})$.

If instead $M_{tt} = 1$ for all $t$, then $\tr(M) = N = \tr(|M|)$, which gives $\displaystyle\sum_{i=1}^N (|\lambda_i| - \lambda_i) = 0$. Since $|\lambda_i| \ge \lambda_i$ for all real $\lambda_i$, every eigenvalue must be non-negative ($\lambda_i \ge 0$), meaning $M$ is positive semidefinite. For any three vertices with $M_{ij} = 1$ and $M_{jk} = 1$, the determinant of their $3 \times 3$ principal submatrix is $-(M_{ik} - 1)^2 \ge 0$, which forces $M_{ik} = 1$. Because $M$ is connected, this transitivity forces all entries to be $1$, uniquely establishing $M \cong J_N$.
\end{proof}

\begin{corollary}\label{cor:kronecker_invariant_classification}
For any simple graph $G$, the operation $B \mapsto M \kron \A(G)$ preserves the vertex energy of $G$ at every index ($\E_B((t, i)) = \E_G(v_i)$) if and only if
\[
M \cong \left( \bigoplus_{i=1}^{c_1} \A(K_{r_i, r_i}) \right) \oplus \left( \bigoplus_{j=1}^{c_2} J_{k_j} \right) \qquad (r_i, k_j \ge 1, \; c_1 + c_2 \ge 1).
\]
\end{corollary}
\begin{proof}
By Theorem~\ref{thm:master}, $\E_B((t, i)) = \E_G(v_i)$ holds for all $t, i$ if and only if $|M|_{tt} = 1$ for all $t$. Decomposing $M \cong \displaystyle\bigoplus_{\ell=1}^c M_\ell$ into its connected components gives $|M| \cong \displaystyle\bigoplus_{\ell=1}^c |M_\ell|$, which satisfies $|M|_{tt} = 1$ if and only if $|M_\ell|_{ss} = 1$ for every slot $s$ of each component. Since $M$ has no isolated zero-diagonal slots, applying Theorem~\ref{thm:connected_matrix_classification} to each connected component implies $M_\ell \cong \A(K_{r_i, r_i})$ ($r_i \ge 1$) or $M_\ell \cong J_{k_j}$ ($k_j \ge 1$). Grouping these blocks diagonally completes the proof.
\end{proof}

\section{Conclusion}
We used the Kronecker product framework to provide a simple and unified way to understand how vertex energy changes under graph operations and when it remains invariant. For future work, it would be interesting to see whether vertex energy can help distinguish non-isomorphic cospectral graphs. Another natural problem is to find operations that produce graphs where every vertex has the same energy, and to study how vertex energy behaves under other graph operations.

\end{document}